\documentclass[12pt]{amsart}
\usepackage[T1]{fontenc}
\usepackage[utf8]{inputenc}
\usepackage{lmodern}
\usepackage[margin=1.05in]{geometry}
\usepackage{amsmath,amssymb,amsthm,mathtools}
\usepackage{microtype}
\usepackage{enumitem}
\usepackage{xcolor}
\usepackage{tikz}
\usepackage[colorlinks=true,linkcolor=blue!50!black,citecolor=blue!50!black,urlcolor=blue!50!black]{hyperref}

\newtheorem{theorem}{Theorem}[section]
\newtheorem{proposition}[theorem]{Proposition}
\newtheorem{corollary}[theorem]{Corollary}
\newtheorem{lemma}[theorem]{Lemma}
\theoremstyle{definition}
\newtheorem{definition}[theorem]{Definition}
\theoremstyle{remark}

\DeclareMathOperator{\reg}{reg}
\DeclareMathOperator{\creg}{creg}
\DeclareMathOperator{\Ind}{Ind}

\newcommand{\N}{\mathbb N}
\newcommand{\kk}{\mathrm{k}}
\newcommand{\Gens}{\mathcal G}
\newcommand{\Htilde}{\widetilde H}
\newcommand{\bd}{\mathbf d}
\newcommand{\bc}{\mathbf c}
\newcommand{\bq}{\mathbf q}

\setlist[enumerate]{label=\textup{(\roman*)},leftmargin=2.2em}
\title[Nonlinear syzygies of chordal cover powers]{Nonlinear syzygies in ordinary powers\\of chordal cover ideals}

\author[N.T. Hang]{Nguyen Thu Hang}
\address{Thai Nguyen University of Sciences, Phan Dinh phung Ward, Thai Nguyen, Vietnam}
\email{hangnt@tnus.edu.vn}

\author{Thanh Vu}
\address{Institute of Mathematics, VAST, 18 Hoang Quoc Viet, Hanoi, Vietnam}
\email{vuqthanh@gmail.com}


\begin{document}

\begin{abstract}
We construct connected chordal graphs whose cover ideals have ordinary powers that are not componentwise linear over any field, disproving the Herzog--Hibi--Ohsugi conjecture.
\end{abstract}

\maketitle

\section{Introduction}
A homogeneous ideal $I$ in a polynomial ring $S$ over a field $\mathrm{k}$ is componentwise linear if the ideal $I_{\langle d \rangle}$ generated by the forms of degree $d$ in $I$ has a linear resolution for every $d$. This notion was introduced by Herzog and Hibi \cite{HH}. A squarefree monomial ideal is componentwise linear if and only if its Alexander dual is sequentially Cohen--Macaulay. Francisco and Van Tuyl \cite{FV} proved that edge ideals of chordal graphs are sequentially Cohen--Macaulay, which is equivalent to the statement that their cover ideals are componentwise linear. Subsequently, Herzog, Hibi, and Ohsugi \cite[Conjecture~2.5]{HHO} conjectured that all powers of the cover ideals of chordal graphs are componentwise linear, and they verified this conjecture for Cohen--Macaulay chordal graphs. Later, Kumar and Kumar proved the assertion for trees \cite{KK}, and Lu and Wang established the componentwise linearity of all ordinary cover powers for fully clique-whiskered graphs \cite{LW}. For further results, we refer the reader to the survey article by Ha and Van Tuyl \cite{HaVT}. Corresponding problems have also been considered for symbolic powers; for partial results in this direction, see Selvaraja and Skelton \cite{SS}, while an example of a chordal graph whose symbolic powers of cover ideals are not componentwise linear is provided in \cite[Remark~5.12]{DHNT}.

Since powers are generally not squarefree, their homological behavior is no longer directly governed by the independence complex. This raises a natural structural question: can a chordal graph conceal combinatorial information that is invisible to the nonlinearity of its cover ideals yet becomes detectable in higher powers? Our construction relies on this philosophy. Specifically, we associate a chordal graph $H(F)$ to an auxiliary graph $F$ that transfers topological obstructions from the independence complex of $F$ to the multigraded strands of the resolutions of powers of $J(H(F))$.

Our construction is as follows. Let $F$ be a graph on $n$ vertices whose edge set $E(F)$ consists of $m$ edges. For each edge $e \in E(F)$, we introduce two new vertices $y_e$ and $z_e$. We define $H(F)$ to be the graph on the vertex set
$$V(H(F)) = \{x_1,\ldots,x_n\} \cup \{y_e,z_e \mid e\in E(F)\}$$
where $x_1,\ldots,x_n$ form a clique, and for each edge $e=\{i,j\}$ of $F$, we add the edges $\{x_i,y_e\}$, $\{x_j,y_e\}$, and $\{y_e,z_e\}$. The resulting graph $H(F)$ is connected and chordal for any graph $F$. Each part of this construction plays a specific role. Every minimum vertex cover of $H(F)$ has size $d = n + m - 1$; such a cover omits exactly one clique vertex and chooses exactly one vertex from each pair $\{y_e,z_e\}$. Omitting $x_i$ forces the inclusion of $y_e$ for every edge $e$ incident to $i$. Thus, a minimum vertex cover records a single vertex of $F$, whereas a product of minimum vertex cover monomials records a collection of vertices with multiplicities. Carefully chosen multidegrees impose compatibility conditions on these omissions. In this way, the independent sets of $F$ emerge naturally from the ordinary powers of the cover ideal.

The following is the fundamental transfer lemma that allows us to control the syzygies of a graded component of the cover ideal of $H(F)$ via the combinatorial information of $F$. Let 
\[
\Delta_s(F) = \{ A \subseteq [n] \mid A \text{ is independent in } F \text{ and } |A| = s \}
\]
be the pure simplicial complex generated by the independent sets of $F$ of cardinality $s$. Suppose that $F$ has no isolated vertices and that $\Delta_s(F)$ is non-empty. Set 
\[
L_s = (J(H(F))^s)_{\langle sd \rangle} \quad \text{and} \quad M_s = \prod_{i=1}^n x_i^s \prod_{e \in E(F)} y_e z_e^{s-1}.
\]
We show that the upper Koszul complex of $L_s$ in multidegree $M_s$ is precisely $\Delta_s(F)$. Consequently, 
$$\beta_{i,M_s}(L_s) = \dim_{\mathrm{k}} \widetilde{H}_{i-1}\bigl(\Delta_s(F); \mathrm{k}\bigr).$$
This correspondence also recovers the homology of links by replacing $M_s$ with suitable monomial factors. This identity explains precisely which topological features obstruct componentwise linearity. Since $\deg M_s = sd + s$, the homology of $\Delta_s(F)$ below dimension $s-1$ produces Betti numbers outside the linear strand. Applying the same observation to links yields the following:

\begin{theorem}\label{thm_CM}
If $\Delta_s(F)$ is not Cohen-Macaulay over $\mathrm{k}$, then $J(H(F))^r$ is not componentwise linear for all $r \ge s$.
\end{theorem}

\begin{figure}[htbp]
\centering
\begin{tikzpicture}[scale=.67,
 core/.style={circle,draw,fill=blue!8,minimum size=8mm,inner sep=1pt},
 aux/.style={circle,draw,fill=white,minimum size=8mm,inner sep=1pt},
 every node/.style={font=\small}]
\node[core] (x1) at (-1.5,1.5) {$x_1$};
\node[core] (x2) at (1.5,1.5) {$x_2$};
\node[core] (x3) at (1.5,-1.5) {$x_3$};
\node[core] (x4) at (-1.5,-1.5) {$x_4$};
\draw (x1)--(x2)--(x3)--(x4)--(x1);
\draw[gray] (x1)--(x3) (x2)--(x4);
\node[aux] (y12) at (0,3.3) {$y_{12}$};
\node[aux] (z12) at (0,4.8) {$z_{12}$};
\node[aux] (y23) at (3.3,0) {$y_{23}$};
\node[aux] (z23) at (4.8,0) {$z_{23}$};
\node[aux] (y34) at (0,-3.3) {$y_{34}$};
\node[aux] (z34) at (0,-4.8) {$z_{34}$};
\node[aux] (y41) at (-3.3,0) {$y_{41}$};
\node[aux] (z41) at (-4.8,0) {$z_{41}$};
\draw (x1)--(y12)--(x2) (y12)--(z12);
\draw (x2)--(y23)--(x3) (y23)--(z23);
\draw (x3)--(y34)--(x4) (y34)--(z34);
\draw (x4)--(y41)--(x1) (y41)--(z41);
\end{tikzpicture}
\caption{The connected chordal graph $H(K_{2,2})$.}
\label{fig:graph}
\end{figure}
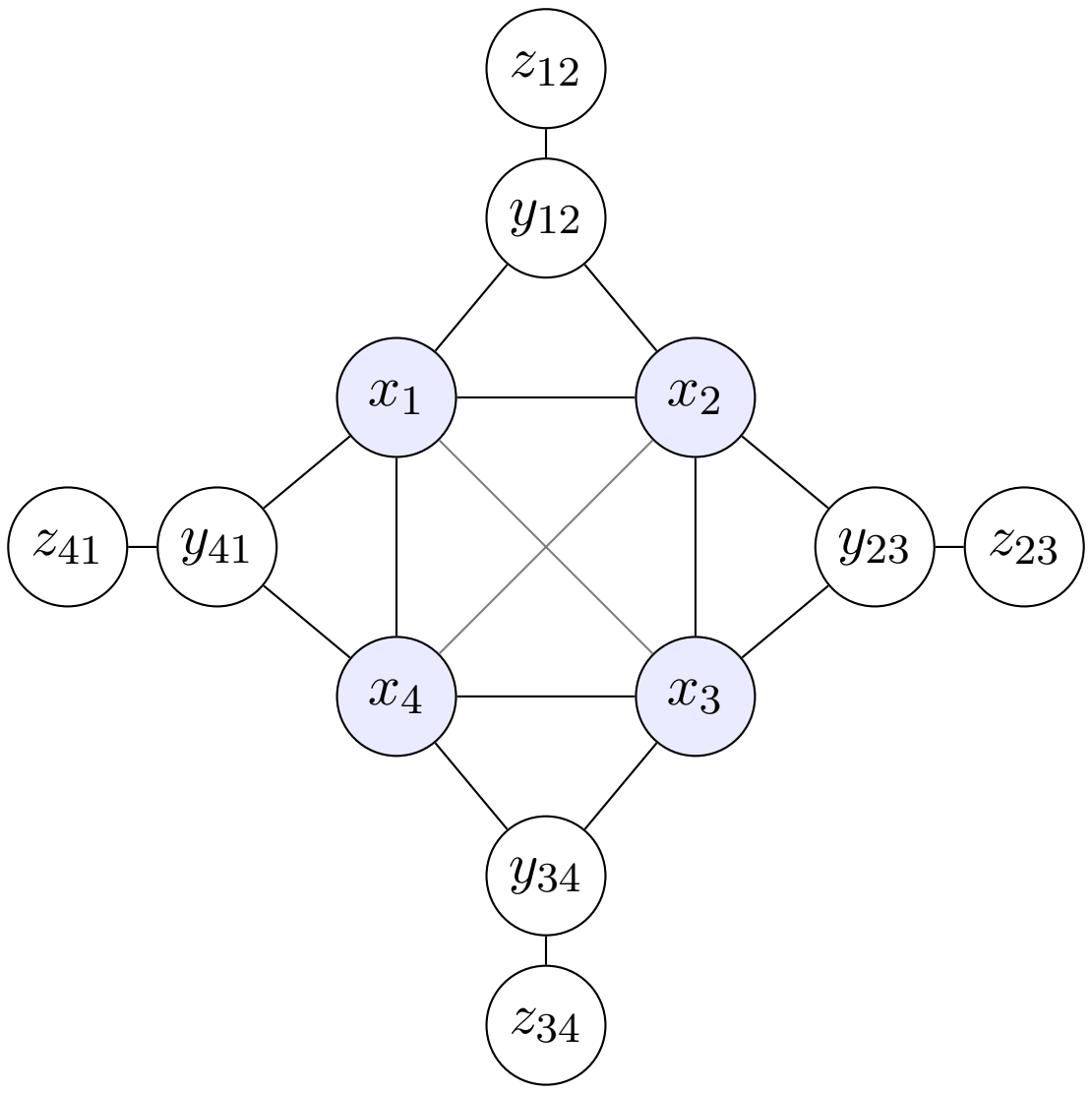

The simplest example arises from $F = K_{2,2}$, for which $\Delta_2(F)$ has non-vanishing reduced zeroth homology. This implies that the graph $H(K_{2,2})$ has $12$ vertices and that $J(H(K_{2,2}))^r$ is not componentwise linear for all $r \ge 2$. This conclusion holds over every field. Figure~\ref{fig:graph} illustrates the graph $H(K_{2,2})$.

The significance of this construction lies in identifying the precise mechanism responsible for the failure of componentwise linearity for powers of cover ideals of chordal graphs.

Section~\ref{sec:prelim} collects the necessary notation and homological tools. Section~\ref{sec:construction} describes the cover generators and the least-degree components. We then establish the transfer result for homology. Finally, the Cohen--Macaulay consequences and further applications are derived in Section~\ref{sec:CM}.

\section{Preliminaries}\label{sec:prelim}

Throughout this paper, let $\kk$ be a field, and let $[n] = \{1, \ldots, n\}$. For a vector $\mathbf{a} \in \mathbb{N}^n$, we write $|\mathbf{a}| = \sum_i a_i$. For a subset $B \subseteq [n]$, let $\mathbf{1}_B$ denote its characteristic vector, and let $x_B = \prod_{i \in B} x_i$.

\subsection{Graphs and cover ideals}

All graphs considered are finite and simple. For a graph $G$, let $V(G)$ and $E(G)$ denote its vertex set and edge set, respectively. A subset of vertices is a \emph{clique} if its vertices are pairwise adjacent, and it is \emph{independent} if no two of its vertices are adjacent. The independence number $\alpha(G)$ is the maximum cardinality of an independent set in $G$. For $A \subseteq V(G)$, its closed neighborhood is defined as
\[
N_G[A] = A \cup \{v \in V(G) \mid \{u, v\} \in E(G) \text{ for some } u \in A\},
\]
and $G \setminus A$ denotes the induced subgraph on $V(G) \setminus A$.

A graph $G$ is \emph{chordal} if it contains no induced cycles of length four or greater. A vertex of $G$ is \emph{simplicial} if its neighbors form a clique.

A \emph{vertex cover} of $G$ is a subset of vertices that intersects every edge of $G$. A vertex cover is \emph{minimal} if no proper subset is a vertex cover. In the polynomial ring $S = \kk[w \mid w \in V(G)]$, the \emph{cover ideal} $J(G)$ is defined as
\[
J(G) = \bigcap_{\{u,v\} \in E(G)} (u, v) = \left( \prod_{w \in C} w \;\Bigg|\; C \text{ is a minimal vertex cover of } G \right).
\]

\subsection{Betti numbers and componentwise linearity}
For a nonzero homogeneous ideal $I$, let
\[
\alpha(I) = \min\{d \mid I_d \neq 0\}.
\]
For each $d \ge 0$, let $I_{\langle d \rangle} = (I_d)$ denote the ideal generated by the forms of degree $d$ in $I$. The minimal monomial generating set of a monomial ideal $L$ is denoted by $\Gens(L)$. 

For a monomial ideal $L \subseteq S$ and a monomial $M$, the multigraded Betti number of $L$ in degree $M$ is defined by
\[
\beta_{p,M}(L) = \dim_{\mathrm{k}} \operatorname{Tor}_p^S(L, \mathrm{k})_M.
\]
The total graded Betti numbers and the Castelnuovo--Mumford regularity are given by
\[
\beta_{p,j}(L) = \sum_{\deg M = j} \beta_{p,M}(L), \qquad \reg(L) = \max\{j - p \mid \beta_{p,j}(L) \neq 0\}.
\]

An ideal generated in a single degree $a$ has a \emph{$d$-linear resolution} if $\beta_{p,j}(L) = 0$ for $j \neq a + p$; equivalently, $\reg(L) = a$.

A homogeneous ideal $I$ is componentwise linear if every nonzero component $I_{\langle d \rangle}$ has a $d$-linear resolution \cite{HH}. We use the notion of the componentwise regularity introduced by Caviglia and Varbaro \cite[Definition~1.1]{CV}:

\begin{definition}
Let $I$ be a homogeneous $S$-ideal. The \emph{componentwise regularity} of $I$ is defined as
\[
\operatorname{creg}(I) = \sup_{d : I_d \neq 0} \{\reg(I_{\langle d \rangle}) - d\}.
\]
\end{definition}

In particular, $I$ is componentwise linear if and only if $\operatorname{creg}(I) = 0$.

\subsection{Simplicial complexes}

A simplicial complex $\Delta$ is a collection of subsets of a finite vertex set that is closed under inclusion. Its facets are its maximal faces, and its dimension is defined as $\dim \Delta = \max\{|B| - 1 \mid B \in \Delta\}$. A complex is pure if all its facets have the same cardinality. The pure $(k-1)$-skeleton $\Delta^{[k-1]}$ is the subcomplex generated by all faces of cardinality $k$; whenever such faces exist, it is pure of dimension $k-1$. 

For a face $A \in \Delta$, its link is defined by
\[
\operatorname{link}_\Delta(A) = \{B \mid B \cap A = \emptyset, \ A \cup B \in \Delta\}.
\]

We use reduced simplicial homology $\widetilde{H}_j(\Delta; \mathrm{k})$ with coefficients in the field $\mathrm{k}$. The void complex has a trivial augmented chain complex, whereas $\widetilde{H}_{-1}(\{\varnothing\}; \mathrm{k}) \cong \mathrm{k}$.

For a simplicial complex $\Delta$ on the vertex set $[N]$, its Stanley--Reisner ring is given by
\[
\mathrm{k}[\Delta] = \mathrm{k}[w_1, \ldots, w_N] / I_\Delta,
\qquad
I_\Delta = \left(\prod_{i \in B} w_i \;\Bigg|\; B \subseteq [N], \ B \notin \Delta\right).
\]
The complex $\Delta$ is said to be Cohen--Macaulay over $\mathrm{k}$ if its Stanley--Reisner ring $\mathrm{k}[\Delta]$ is Cohen--Macaulay, meaning that its depth equals its Krull dimension.

\subsection{Upper Koszul complexes}

For a monomial ideal $L\subseteq\kk[w_1,\ldots,w_N]$ and a monomial $M$,
define
\[
\mathcal K^M(L)=
\left\{\sigma\subseteq\{1,\ldots,N\}:
 \prod_{i\in\sigma}w_i\mid M,\
 \frac{M}{\prod_{i\in\sigma}w_i}\in L\right\}.
\]
This is the upper Koszul complex of $L$ at $M$. The following is \cite[Theorem~2.2]{BCP}.

\begin{lemma}\label{lem:koszul}
For every $p\ge0$,
\[
\beta_{p,M}(L)=
\dim_{\kk}\Htilde_{p-1}(\mathcal K^M(L);\kk).
\]
\end{lemma}

\section{The chordal construction}\label{sec:construction}

Let $F$ be a graph on $[n] = \{1, \ldots, n\}$, with $n \ge 2$ and no isolated vertices. Let $m = |E(F)|$.

\begin{definition}\label{def:construction}
The graph $H(F)$ has the vertex set
\[
V(H(F)) = \{x_1, \ldots, x_n\} \cup \{y_e, z_e \mid e \in E(F)\}.
\]
The vertices $x_1, \ldots, x_n$ form a complete graph (clique). For every edge $e = \{i, j\} \in E(F)$, the additional edges are
\[
\{x_i, y_e\}, \quad \{x_j, y_e\}, \quad \text{and} \quad \{y_e, z_e\}.
\]
There are no other edges in $H(F)$.
\end{definition}

\begin{proposition}\label{prop:chordal}
The graph $H(F)$ is connected, chordal, and has $n + 2m$ vertices.
\end{proposition}

\begin{proof}
The graph $H(F)$ is connected because each $z_e$ is connected to $y_e$, each $y_e$ is connected to $x_i$ and $x_j$ (which belong to the central clique), and $F$ has no isolated vertices. 

To show that $H(F)$ is chordal, let $C$ be an induced cycle of length at least $4$ in $H(F)$. Since each $z_e$ is a leaf in $H(F)$, no vertex $z_e$ can belong to $C$. If $C$ is entirely contained within the clique $\{x_1, \ldots, x_n\}$, its length can be at most $3$, which contradicts our assumption. Hence, $C$ must contain at least one vertex $y_e$ from the set $\{y_e \mid e \in E(F)\}$ for some edge $e = \{i, j\}$. The only neighbors of $y_e$ in $H(F) \setminus \{z_e\}$ are $x_i$ and $x_j$; consequently, $C$ must contain both $x_i$ and $x_j$. But then $y_e, x_i, x_j$ form a triangle (since $x_i$ and $x_j$ are connected by an edge in the clique), which forces $C$ to be this $3$-cycle. This contradicts the assumption that $C$ has length at least $4$. Therefore, no induced cycle of length $4$ or greater can exist, proving that $H(F)$ is chordal.
\end{proof}

Let
\[
S_F = \mathrm{k}[x_1, \ldots, x_n, y_e, z_e \mid e \in E(F)],
\quad I_F = J(H(F)), \quad D = n + m - 1.
\]

\begin{proposition}\label{prop:covers}
The minimal vertex covers of $H(F)$ are precisely the following two types:
\begin{enumerate}
\item Choose $i \in [n]$ and omit $x_i$, while including all other core vertices $x_j$ ($j \neq i$). For each edge $e$ incident to $i$, include $y_e$ and omit $z_e$. For every remaining edge $e$, choose exactly one of $y_e$ or $z_e$. These covers have cardinality $D$.
\item Include all core vertices, and choose exactly one of $y_e$ or $z_e$ for each edge $e \in E(F)$. The set of edges $e$ for which $z_e$ is chosen must form an edge cover of $F$ (i.e., it must meet every vertex of $F$). These covers have cardinality $D + 1$.
\end{enumerate}
Consequently, the minimal generators of $I_F$ have degrees $D$ and $D + 1$, and $\alpha(I_F) = D$.
\end{proposition}

\begin{proof}
Any vertex cover of $H(F)$ must contain at least $n - 1$ core vertices and at least one vertex from each pair $\{y_e, z_e\}$ for every edge $e \in E(F)$. Moreover, a minimal cover cannot contain both $y_e$ and $z_e$ because $z_e$ is a leaf adjacent only to $y_e$.

We analyze the possibilities based on the number of core vertices included:
\begin{enumerate}
\item[(i)] \emph{Omitting one core vertex:} Suppose $x_i$ is omitted. To cover the edges between $x_i$ and its incident edges $y_e$, we are forced to include $y_e$ (and thus omit $z_e$) for every edge $e$ incident to $i$. For all remaining edges, we may independently choose either $y_e$ or $z_e$. Every included core vertex is necessary to cover the edges incident to the omitted vertex $x_i$, and each chosen vertex from $\{y_e, z_e\}$ is necessary to cover the edge $\{y_e, z_e\}$. Thus, these sets form minimal vertex covers of cardinality $(n - 1) + m = D$.

\item[(ii)] \emph{Including all core vertices:} Suppose all core vertices $x_1, \ldots, x_n$ are included. A core vertex $x_i$ is necessary in the cover if and only if some adjacent vertex $y_e$ is omitted, which forces the choice of $z_e$ for some edge $e$ incident to $i$. Therefore, all core vertices are simultaneously necessary if and only if the set of edges $e$ for which $z_e$ is chosen covers every vertex of $F$ (i.e., forms an edge cover of $F$). Each chosen vertex from the pairs $\{y_e, z_e\}$ is again necessary. This yields minimal vertex covers of cardinality $n + m = D + 1$.
\end{enumerate}

The degree assertions follow directly from the cardinalities of these minimal vertex covers.
\end{proof}

For $s \ge 1$, let
\[
L_s = (I_F^s)_{\langle sD \rangle}.
\]

\begin{theorem}\label{thm:generators}
The minimal monomial generators of $L_s$ are precisely the monomials of the form
\begin{equation}\label{eq:generators}
u(\bc, \bq) = \prod_{i=1}^n x_i^{s - c_i} \prod_{e \in E(F)} y_e^{q_e} z_e^{s - q_e},
\end{equation}
where $\bc \in \mathbb{N}^n$ and $\bq = (q_e)_{e \in E(F)} \in \mathbb{N}^{E(F)}$ satisfy
\begin{equation}\label{eq:constraints}
|\bc| = s \quad \text{and} \quad c_i + c_j \le q_e \le s \quad \text{for each edge } e = \{i, j\} \in E(F).
\end{equation}
Moreover, distinct admissible pairs $(\bc, \bq)$ yield distinct monomials.
\end{theorem}

\begin{proof}
Every monomial of degree $sD$ in $I_F^s$ must be a product of $s$ generators of $I_F$, each of degree $D$. By Proposition~\ref{prop:covers}, the degree-$D$ generators of $I_F$ are precisely the covers of type (i), each of which omits exactly one core vertex $x_i$. When we multiply $s$ such generators together, let $c_i$ denote the number of factors that omit the vertex $x_i$. Then we must have $|\bc| = s$, and the exponent of $x_i$ in the resulting product is $s - c_i$. 

For a given edge $e = \{i, j\}$, a factor must include $y_e$ if and only if it omits either $x_i$ or $x_j$. Thus, exactly $c_i + c_j$ of the $s$ factors are forced to choose $y_e$. Among the remaining $s - (c_i + c_j)$ factors, we may independently choose either $y_e$ or $z_e$. Consequently, the total exponent $q_e$ of $y_e$ satisfies the constraints in \eqref{eq:constraints}, and the exponent of $z_e$ is necessarily $s - q_e$.

Conversely, let us fix any admissible pair $(\bc, \bq)$. We can construct a product of $s$ generators of type (i) as follows: assign $c_i$ of the $s$ factor positions to omit $x_i$ for each $i \in [n]$. For each edge $e = \{i, j\}$, $y_e$ is automatically included in the $c_i + c_j$ forced positions. For the remaining factors, we choose $y_e$ in exactly $q_e - c_i - c_j$ positions and $z_e$ in the rest. Since these choices can be made independently for each edge, this produces a valid collection of $s$ generators of type (i) whose product is precisely $u(\bc, \bq)$ given in \eqref{eq:generators}.

Finally, the exponents of the core variables uniquely recover $\bc$, and the exponents of the $y_e$ variables uniquely recover $\bq$, ensuring that distinct admissible pairs determine distinct minimal generators.
\end{proof}

We now prove the main transfer result for the Betti numbers of $L_s$ at a specified multidegree.
\begin{theorem}\label{thm:transfer}
Let $F$ be a graph with no isolated vertices, and let $A$ be an independent set of $F$. Let $k \ge 1$, and suppose that $F \setminus N_F[A]$ has an independent set of cardinality $k$. Fix $s \ge |A| + k$ and $\mathbf{d} \in \mathbb{N}^n$ with
\[
|\mathbf{d}| = s - |A| - k.
\]
For $a_i = 1$ if $i \in A$ and $0$ otherwise, define the multidegree
\begin{equation}\label{eq:multidegree}
M_{s,A,k,\mathbf{d}} = \prod_{i=1}^n x_i^{s - d_i - a_i} \prod_{e=\{i,j\} \in E(F)} y_e^{d_i + d_j + 1} z_e^{s - d_i - d_j - 1}.
\end{equation}
Let $\Gamma = \Delta_{k}(F \setminus N_F[A])$ denote the $(k-1)$-skeleton of the independence complex of $F \setminus N_F[A]$. Then, writing $M = M_{s,A,k,\mathbf{d}}$:
\begin{enumerate}
\item $M$ is a monomial of total degree $sD + k$;
\item the minimal generators of $L_s$ dividing $M$ are precisely
\begin{equation}\label{eq:divisors}
\left\{ \frac{M}{x_B} \;\Bigg|\; B \subseteq V(F) \setminus N_F[A], \ B \text{ is independent in } F, \ |B| = k \right\};
\end{equation}
\item the upper Koszul complex $\mathcal{K}^M(L_s)$ is naturally isomorphic to $\Gamma$;
\item for every $p \ge 0$,
\begin{equation}\label{eq:transfer}
\beta_{p,M}(L_s) = \dim_{\mathrm{k}} \widetilde{H}_{p-1}(\Gamma; \mathrm{k}).
\end{equation}
\end{enumerate}
\end{theorem}

\begin{proof}
First, since $d_i + d_j \le |\mathbf{d}| = s - |A| - k$, the exponent of $z_e$ in $M$ is at least $|A| + k - 1 \ge 0$. The exponent of $y_e$ is $d_i + d_j + 1 \ge 1$. If $a_i = 0$, the exponent of $x_i$ is $s - d_i \ge |A| + k \ge 1$; if $a_i = 1$, it is $s - d_i - 1 \ge |A| + k - 1 \ge 1$. Thus, $M$ is a well-defined monomial. Its core-variable degree is
\[
\sum_{i=1}^n (s - d_i - a_i) = ns - (s - |A| - k) - |A| = (n - 1)s + k.
\]
Each pair $y_e, z_e$ contributes a total degree of $s$, so the total degree of $M$ is $(n - 1)s + k + m \cdot s = s(n + m - 1) + k = sD + k$. This proves assertion (i).

To prove assertion (ii), consider a minimal generator $u(\mathbf{c}, \mathbf{q})$ of $L_s$ (from Theorem~\ref{thm:generators}) that divides $M$. For every edge $e = \{i, j\}$, both $u(\mathbf{c}, \mathbf{q})$ and $M$ have total degree $s$ on the variable pair $\{y_e, z_e\}$. Divisibility therefore forces equality of the exponents on each pair, yielding
\begin{equation}\label{eq:qforced}
q_e = d_i + d_j + 1.
\end{equation}
Divisibility in the core variables implies $c_i \ge d_i + a_i$ for all $i$. Setting $b_i = c_i - d_i$, the constraints on the generator $u(\mathbf{c}, \mathbf{q})$ translate to
\begin{equation}\label{eq:bconditions}
b_i \ge a_i, \quad \sum_{i=1}^n b_i = |A| + k, \quad b_i + b_j \le 1 \quad \text{for every edge } \{i, j\} \in E(F).
\end{equation}
Since every vertex of $F$ has at least one incident edge (as $F$ has no isolated vertices) and all $b_i$ are non-negative integers, the inequality $b_i + b_j \le 1$ implies that $b_i \le 1$ for all $i$. Consequently, the support of $\mathbf{b}$ corresponds to an independent set of cardinality $|A| + k$ containing $A$. This set can be written as $A \cup B$, where $B$ is an independent $k$-set in $F \setminus N_F[A]$. For this choice, we have
\[
c_i = d_i + a_i + \mathbf{1}_B(i).
\]
Combining this with \eqref{eq:qforced}, we obtain $u(\mathbf{c}, \mathbf{q}) = M / x_B$.

Conversely, let $B$ be any independent $k$-set in $F \setminus N_F[A]$. Setting $\mathbf{c} = \mathbf{d} + \mathbf{1}_A + \mathbf{1}_B$ and determining $\mathbf{q}$ via \eqref{eq:qforced}, we find that $|\mathbf{c}| = s$. The independence of $A \cup B$ ensures that $c_i + c_j \le d_i + d_j + 1 = q_e$, and $q_e \le |\mathbf{d}| + 1 \le s$. Thus, the pair $(\mathbf{c}, \mathbf{q})$ is admissible, proving assertion (ii).

For assertion (iii), a face $\sigma$ belongs to the upper Koszul complex $\mathcal{K}^M(L_s)$ if and only if $M / \prod_{w \in \sigma} w$ is divisible by some minimal generator of $L_s$ dividing $M$. By assertion (2), these generators are precisely of the form $M / x_B$, where $B$ is an independent $k$-set in $F \setminus N_F[A]$ disjoint from $A$. Thus, divisibility holds if and only if $\sigma$ consists of core variables indexed by a subset of $B$. The faces of the upper Koszul complex are therefore generated by these independent $k$-sets, which is precisely the $(k-1)$-skeleton $\Gamma = \Delta_{k}(F \setminus N_F[A])$. 

Finally, assertion (iv) follows immediately from the identification of $\mathcal{K}^M(L_s)$ with $\Gamma$ and Lemma~\ref{lem:koszul}.
\end{proof}

\section{Cohen--Macaulay obstructions}\label{sec:CM}
We now deduce the Cohen--Macaulay obstructions for the componentwise linearity of powers of the cover ideal of $H(F)$. Recall that $\Delta_q(F)$ denotes the $(q-1)$-skeleton of the independence complex of $F$.

\begin{theorem}\label{thm:CM}
Let $F$ be a graph with no isolated vertices, and let $s \ge 1$. If $I_F^s$ is componentwise linear over $\mathrm{k}$, then
\[
\Delta_q(F) \quad \text{is Cohen--Macaulay over } \mathrm{k} \quad \text{for every } 1 \le q \le \min\{s, \alpha(F)\}.
\]
\end{theorem}

\begin{proof}
Suppose that $\Gamma = \Delta_{q-1}(F)$ is not Cohen--Macaulay over $\mathrm{k}$. Since $\Gamma$ is pure of dimension $q-1$, Reisner's criterion implies that there exist a face $A \in \Gamma$ and an integer $j < \dim \operatorname{link}_\Gamma(A)$ such that
\[
\widetilde{H}_j(\operatorname{link}_\Gamma(A); \mathrm{k}) \neq 0.
\]
Let $k = q - |A|$. Every facet of $\Gamma$ containing $A$ has cardinality $q$, so the link $\operatorname{link}_\Gamma(A)$ has dimension $k-1$.

If $k = 0$, the link is $\{\varnothing\}$, which has no non-vanishing reduced homology in degrees $j < -1$. If $k \ge 1$, the link has at least one vertex, and its reduced homology in degree $-1$ vanishes. Thus, any failure of the Cohen--Macaulay property must occur with $j \ge 0$, $k \ge 2$, and
\[
1 \le p := j + 1 < k.
\]

By definition, a face $\sigma$ belongs to $\operatorname{link}_\Gamma(A)$ if and only if $\sigma$ is an independent set in $F \setminus N_F[A]$ of cardinality at most $k$. Therefore, we have the equality
\[
\operatorname{link}_\Gamma(A) = \Delta_{k}(F \setminus N_F[A]).
\]

Now, fix any $s \ge q = |A| + k$. Choose $\mathbf{d} \in \mathbb{N}^n$ such that $|\mathbf{d}| = s - q = s - |A| - k$, and apply Theorem~\ref{thm:transfer} with this choice of $A$, $k$, and $\mathbf{d}$. The transfer formula \eqref{eq:transfer} yields
\[
\beta_{p, M}(L_s) = \dim_{\mathrm{k}} \widetilde{H}_{p-1}\bigl(\Delta_{k-1}(F \setminus N_F[A]); \mathrm{k}\bigr) = \dim_{\mathrm{k}} \widetilde{H}_{j}\bigl(\operatorname{link}_\Gamma(A); \mathrm{k}\bigr) \neq 0,
\]
where $M = M_{s,A,k,\mathbf{d}}$. By assertion (1) of Theorem~\ref{thm:transfer}, the total degree of this multidegree is
\[
\deg M = sD + k.
\]
Since $p < k$, we have $sD + k > sD + p$, which implies that the graded component $L_s = (I_F^s)_{\langle sD \rangle}$ has a nonzero Betti number outside its linear strand. Consequently, $I_F^s$ is not componentwise linear. This completes the proof.
\end{proof}

\begin{proof}[Proof of Theorem \ref{thm_CM}]
The result follows immediately from Theorem \ref{thm:CM}.
\end{proof}

We conclude the paper with an example showing that $\operatorname{creg}(I_F^s)$ can be arbitrarily large. For $r \ge 2$, let
\[
G_r = H(K_{r,r}), \quad I_r = J(G_r), \quad D_r = r^2 + 2r - 1, \quad L_r = (I_r^s)_{\langle sD_r \rangle}.
\]
The graph $G_r$ is connected, chordal, and has $2r(r + 1)$ vertices.

\begin{theorem}\label{thm:Krr}
Let $r, s \ge 2$. For every integer $k$ with $2 \le k \le \min\{r, s\}$, there exists a monomial $M$ of degree $sD_r + k$ such that
\[
\beta_{1, M}\bigl((I_r^s)_{\langle sD_r \rangle}\bigr) = 1.
\]
Consequently:
\begin{enumerate}
\item $I_r^s$ is not componentwise linear over any field for every $s \ge 2$;
\item $\operatorname{reg}(L_r) \ge sD_r + \min\{r, s\} - 1$;
\item $\operatorname{creg}(I_r^s) \ge \min\{r, s\} - 1$.
\end{enumerate}
\end{theorem}

\begin{proof}
Let $U$ and $V$ denote the two partitions of $K_{r,r}$. For every $2 \le k \le r$, we have
\[
\dim_{\mathrm{k}} \widetilde{H}_0(\Delta_k(K_{r,r}); \mathrm{k}) = 1
\]
over every field. 

Choose any vector $\mathbf{d} \in \mathbb{N}^{2r}$ with $|\mathbf{d}| = s - k$, and take $M = M_{s, \varnothing, k, \mathbf{d}}$. Theorem~\ref{thm:transfer} then yields the asserted Betti number. Taking $k = 2$ proves assertion (i), taking $k = \min\{r, s\}$ proves assertion (ii), and assertion (iii) follows immediately from (ii).
\end{proof}

\noindent\textbf{Acknowledgment.} This work is partially supported by the Vietnam National Foundation for Science and Technology Development under the grant number NAFOSTED/101/10104/2026/1279.

\end{document}